\documentclass{article}
\usepackage{amssymb,amsmath,theorem,euscript}

\input{epsf}

\newcounter{sec}

\def\sm{\smallskip}

\newcounter{punct}[sec]

\def\punct{\refstepcounter{punct}{\arabic{sec}.\arabic{punct}.  }} 

\def\COUNTERS{\addtocounter{sec}{1}
              \setcounter{punct}{0}
          \setcounter{equation}{0}
          \setcounter{theorem}{0}
                  }

\newtheorem{theorem}{Theorem}[sec]
\newtheorem{proposition}[theorem]{Proposition}

\begin{document}

 \def\ov{\overline}
\def\wt{\widetilde}
 \newcommand{\rk}{\mathop {\mathrm {rk}}\nolimits}
\newcommand{\Aut}{\mathop {\mathrm {Aut}}\nolimits}
\newcommand{\Out}{\mathop {\mathrm {Out}}\nolimits}
 \newcommand{\tr}{\mathop {\mathrm {tr}}\nolimits}
  \newcommand{\diag}{\mathop {\mathrm {diag}}\nolimits}
  \newcommand{\supp}{\mathop {\mathrm {supp}}\nolimits}
  \newcommand{\indef}{\mathop {\mathrm {indef}}\nolimits}
  \newcommand{\dom}{\mathop {\mathrm {dom}}\nolimits}
  \newcommand{\im}{\mathop {\mathrm {im}}\nolimits}
 
\renewcommand{\Re}{\mathop {\mathrm {Re}}\nolimits}

\def\Br{\mathrm {Br}}

\def\SL{\mathrm {SL}}
\def\Diag{\mathrm {Diag}}
\def\SU{\mathrm {SU}}
\def\GL{\mathrm {GL}}
\def\U{\mathrm U}
\def\OO{\mathrm O}
 \def\Sp{\mathrm {Sp}}
 \def\SO{\mathrm {SO}}
\def\SOS{\mathrm {SO}^*}
 \def\Diff{\mathrm{Diff}}
 \def\Vect{\mathfrak{Vect}}
\def\PGL{\mathrm {PGL}}
\def\PU{\mathrm {PU}}
\def\PSL{\mathrm {PSL}}
\def\Symp{\mathrm{Symp}}
\def\End{\mathrm{End}}
\def\Mor{\mathrm{Mor}}
\def\Aut{\mathrm{Aut}}
 \def\PB{\mathrm{PB}}
 \def\cA{\mathcal A}
\def\cB{\mathcal B}
\def\cC{\mathcal C}
\def\cD{\mathcal D}
\def\cE{\mathcal E}
\def\cF{\mathcal F}
\def\cG{\mathcal G}
\def\cH{\mathcal H}
\def\cJ{\mathcal J}
\def\cI{\mathcal I}
\def\cK{\mathcal K}
 \def\cL{\mathcal L}
\def\cM{\mathcal M}
\def\cN{\mathcal N}
 \def\cO{\mathcal O}
\def\cP{\mathcal P}
\def\cQ{\mathcal Q}
\def\cR{\mathcal R}
\def\cS{\mathcal S}
\def\cT{\mathcal T}
\def\cU{\mathcal U}
\def\cV{\mathcal V}
 \def\cW{\mathcal W}
\def\cX{\mathcal X}
 \def\cY{\mathcal Y}
 \def\cZ{\mathcal Z}
\def\0{{\ov 0}}
 \def\1{{\ov 1}}
 \def\frA{\mathfrak A}
 \def\frB{\mathfrak B}
\def\frC{\mathfrak C}
\def\frD{\mathfrak D}
\def\frE{\mathfrak E}
\def\frF{\mathfrak F}
\def\frG{\mathfrak G}
\def\frH{\mathfrak H}
\def\frI{\mathfrak I}
 \def\frJ{\mathfrak J}
 \def\frK{\mathfrak K}
 \def\frL{\mathfrak L}
\def\frM{\mathfrak M}
 \def\frN{\mathfrak N} \def\frO{\mathfrak O} \def\frP{\mathfrak P} \def\frQ{\mathfrak Q} \def\frR{\mathfrak R}
 \def\frS{\mathfrak S} \def\frT{\mathfrak T} \def\frU{\mathfrak U} \def\frV{\mathfrak V} \def\frW{\mathfrak W}
 \def\frX{\mathfrak X} \def\frY{\mathfrak Y} \def\frZ{\mathfrak Z} \def\fra{\mathfrak a} \def\frb{\mathfrak b}
 \def\frc{\mathfrak c} \def\frd{\mathfrak d} \def\fre{\mathfrak e} \def\frf{\mathfrak f} \def\frg{\mathfrak g}
 \def\frh{\mathfrak h} \def\fri{\mathfrak i} \def\frj{\mathfrak j} \def\frk{\mathfrak k} \def\frl{\mathfrak l}
 \def\frm{\mathfrak m} \def\frn{\mathfrak n} \def\fro{\mathfrak o} \def\frp{\mathfrak p} \def\frq{\mathfrak q}
 \def\frr{\mathfrak r} \def\frs{\mathfrak s} \def\frt{\mathfrak t} \def\fru{\mathfrak u} \def\frv{\mathfrak v}
 \def\frw{\mathfrak w} \def\frx{\mathfrak x} \def\fry{\mathfrak y} \def\frz{\mathfrak z} \def\frsp{\mathfrak{sp}}
 \def\bfa{\mathbf a} \def\bfb{\mathbf b} \def\bfc{\mathbf c} \def\bfd{\mathbf d} \def\bfe{\mathbf e} \def\bff{\mathbf f}
 \def\bfg{\mathbf g} \def\bfh{\mathbf h} \def\bfi{\mathbf i} \def\bfj{\mathbf j} \def\bfk{\mathbf k} \def\bfl{\mathbf l}
 \def\bfm{\mathbf m} \def\bfn{\mathbf n} \def\bfo{\mathbf o} \def\bfp{\mathbf p} \def\bfq{\mathbf q} \def\bfr{\mathbf r}
 \def\bfs{\mathbf s} \def\bft{\mathbf t} \def\bfu{\mathbf u} \def\bfv{\mathbf v} \def\bfw{\mathbf w} \def\bfx{\mathbf x}
 \def\bfy{\mathbf y} \def\bfz{\mathbf z} \def\bfA{\mathbf A} \def\bfB{\mathbf B} \def\bfC{\mathbf C} \def\bfD{\mathbf D}
 \def\bfE{\mathbf E} \def\bfF{\mathbf F} \def\bfG{\mathbf G} \def\bfH{\mathbf H} \def\bfI{\mathbf I} \def\bfJ{\mathbf J}
 \def\bfK{\mathbf K} \def\bfL{\mathbf L} \def\bfM{\mathbf M} \def\bfN{\mathbf N} \def\bfO{\mathbf O} \def\bfP{\mathbf P}
 \def\bfQ{\mathbf Q} \def\bfR{\mathbf R} \def\bfS{\mathbf S} \def\bfT{\mathbf T} \def\bfU{\mathbf U} \def\bfV{\mathbf V}
 \def\bfW{\mathbf W} \def\bfX{\mathbf X} \def\bfY{\mathbf Y} \def\bfZ{\mathbf Z} \def\bfw{\mathbf w}
 \def\R {{\mathbb R }} \def\C {{\mathbb C }} \def\Z{{\mathbb Z}} \def\H{{\mathbb H}} \def\K{{\mathbb K}}
 \def\N{{\mathbb N}} \def\Q{{\mathbb Q}} \def\A{{\mathbb A}} \def\T{\mathbb T} \def\P{\mathbb P} \def\G{\mathbb G}
 \def\bbA{\mathbb A} \def\bbB{\mathbb B} \def\bbD{\mathbb D} \def\bbE{\mathbb E} \def\bbF{\mathbb F} \def\bbG{\mathbb G}
 \def\bbI{\mathbb I} \def\bbJ{\mathbb J} \def\bbK{\mathbb K} \def\bbL{\mathbb L} \def\bbM{\mathbb M} \def\bbN{\mathbb N} \def\bbO{\mathbb O}
 \def\bbP{\mathbb P} \def\bbQ{\mathbb Q} \def\bbS{\mathbb S} \def\bbT{\mathbb T} \def\bbU{\mathbb U} \def\bbV{\mathbb V}
 \def\bbW{\mathbb W} \def\bbX{\mathbb X} \def\bbY{\mathbb Y} \def\kappa{\varkappa} \def\epsilon{\varepsilon}
 \def\phi{\varphi} \def\le{\leqslant} \def\ge{\geqslant}

\def\UU{\bbU}
\def\Mat{\mathrm{Mat}}
\def\tto{\rightrightarrows}

\def\Gr{\mathrm{Gr}}

\def\graph{\mathrm{graph}}

\def\O{\mathrm{O}}

\def\la{\langle}
\def\ra{\rangle}

\def\B{\mathrm B}
\def\Int{\mathrm{Int}}
\def\LGr{\mathrm{LGr}}


\def\I{\mathbb I}
\def\M{\mathbb M}
\def\T{\mathbb T}
\def\S{\mathrm S}

\def\Lat{\mathrm{Lat}}
\def\LLat{\mathrm{LLat}} 
\def\Mod{\mathrm{Mod}}
\def\LMod{\mathrm{LMod}}
\def\Naz{\mathrm{Naz}}
\def\naz{\mathrm{naz}}
\def\bNaz{\mathbf{Naz}}
\def\AMod{\mathrm{AMod}}
\def\ALat{\mathrm{ALat}}
\def\MAT{\mathrm{MAT}}

\def\Ver{\mathrm{Vert}}
\def\Bd{\mathrm{Bd}}
\def\We{\mathrm{We}}
\def\Heis{\mathrm{Heis}}

\def\bbot{{\bot\!\!\!\bot}}

\begin{center}
\Large\bf

Multiplication of conjugacy classes, colligations, and characteristic functions
of matrix argument

\bigskip

\large\sc
Yury A. Neretin\footnote{Supported by the grants FWF, P22122, P25142.}
\end{center}

{\small We extend the classical construction of operator
colligations and characteristic functions.
Consider the group $G$ of finitary block unitary matrices
of size $\alpha+\infty+\dots+\infty$ ($m$ times). Consider 
the subgroup $K=\U(\infty)$, which consists  of block diagonal unitary matrices
with a block 1 of size $\alpha$ and a matrix $u\in \U(\infty)$ repeated $m$ times.
It appears that there is a natural multiplication
on the  conjugacy classes $G//K$. We construct 'spectral data' of conjugacy classes, 
which visualize the multiplication and are  sufficient for a reconstruction of a conjugacy class.

MSC 22E66. 47A48, 15A72, 16W22
}

\section{Formulation of the problem}

\COUNTERS

{\bf\punct Notation.%
\label{ss:notation}} By $1=1_n$ we denote the unit matrix of order $n$,
by $A^t$ the transposed matrix.
 Let $\Mat(n)$ be the space of complex matrices of size $n\times n$, 
$\GL(n,\C)$ be the group of invertible matrices of order $n$,
$\U(n)$ be the unitary subgroup in $\GL(n,\C)$. 

Denote by $\MAT(\infty)$ the space of all infinite matrices $g$
 such that $g-1$
has only finite number of non-zero elements.
By $\GL(\infty,\C)$, we denote the group of invertible elements of
$\MAT(\infty)$. By $\U(\infty)$
we denote its subgroup consisting of unitary matrices.

Denote $\Mat(\infty)$ the space of all infinite matrices having
only finite number of non-zero matrix elements,
i.e. $g\in\Mat(\infty)\Leftrightarrow g+1\in\MAT(\infty)$.

Let $G$ be a group, $H\subset G$ a subgroup. By $G//H$ we denote the space of conjugacy classes $G$ with respect
to $H$, by $H\setminus G/H$ the double cosets. 

\sm

{\bf\punct  Colligations.%
\label{ss:classic}} The following objects arose  by independent reasons
in spectral theory of non-selfadjoint operators, system theory, and representation theory
of infinite-dimensional groups (see \cite{Liv1}, \cite{Liv2}, \cite{FSN}, \cite{Bro1}, \cite{Bro2}, \cite{Dym},
\cite{Goh},
\cite{Olsh-GB}, \cite{Olsh-CR}, \cite{Ner-book}).

Fix a non-negative integer $\alpha$. 
Denote by $G$ the group $\U(\infty)$, by $K$ its subgroup consisting of matrices of the form
$
\begin{pmatrix}
 1_\alpha&0\\0& h
\end{pmatrix}
$.
Consider the space $G//K$ of conjugacy classes of $G$ with respect to $K$,
i.e., elements of $G$ defined up to the equivalence
$$
\begin{pmatrix}
 a&b\\c&d
\end{pmatrix}
\sim
\begin{pmatrix}
 1_\alpha&0\\0& h
\end{pmatrix}
\begin{pmatrix}
 a&b\\c&d
\end{pmatrix}
\begin{pmatrix}
 1_\alpha&0\\0& h
\end{pmatrix}^{-1}
.$$
There is a natural multiplication on the set $G//K$ defined in the following way. Let $\frg$, $\frh$
be conjugacy classes, let $g$, $h$ be their representatives, let  $g$, $h$ actually be contained in 
$\U(\alpha+N)$, i.e., they are given by block $(\alpha+N+\infty)$-matrices of the form
$$
g=\begin{pmatrix}
 a_1&b_1&0\\c_1&d_1&0\\ 0&0&1_\infty
\end{pmatrix},\qquad g=\begin{pmatrix}
 a_2&b_2&0\\c_2&d_2&0\\ 0&0&1_\infty
\end{pmatrix}
.
$$
We define  the product $\frg\circ\frh$ as the conjugacy class containing
\begin{multline}
\begin{pmatrix}
 a_1&b_1&0\\c_1&d_1&0\\ 0&0&1_\infty
\end{pmatrix}\circ\begin{pmatrix}
 a_2&b_2&0\\c_2&d_2&0\\ 0&0&1_\infty
\end{pmatrix}
:=\\:=
\begin{pmatrix}
 a_1&b_1&0&0\\c_1&d_1&0&0\\0&0&1_N&0\\ 0&0&0&1_\infty
\end{pmatrix}
\begin{pmatrix}
 a_2&0&b_2&0\\0&1_N&0&0\\c_2&0&d_2&0\\0& 0&0&1_\infty
 \end{pmatrix}
 =
 \begin{pmatrix}
 a_1 a_2&b_1&a_1b_2&0\\
 c_1a_2&d_1&c_1 b_2&0\\
 c_2&0&d_2&0\\
 0&0&0&1_\infty
\end{pmatrix}.
\label{eq:classic}
\end{multline}
This product is a well-defined associative operation on $G//K$.

\sm

This operation admits the following visualization:
To a matrix 
$g=\begin{pmatrix}
 a&b\\c&d
\end{pmatrix}$
we assign the {\it characteristic function} as a function on $\C$ taking values in $\alpha\times \alpha$-matrices
given by
$$
\chi(z)=a+z b(1-zd)^{-1} c.
$$
This function depends only on a conjugacy class $\frg$ and not on a matrix $g$, so we can write
$\chi_\frg(z)$. 

Next, denote by $\eta_\frg$ the set of eigenvalues of $d$ that are contained in the circle $|z|=1$.

The following statement is well-known (see, e.g., \cite{Bro2},
in this work
the complete unitary group instead of $\U(\infty)$
is considered).

\begin{theorem}
\label{th:1.1}
{\rm a)} 
$
\chi_{\frg\circ\frh}(z)=\chi_\frg(z)\chi_\frh(z).
$

\sm

{\rm b)} If $|z|\le 1$, then $\|\chi(z)\|\le 1$; if $|z|=1$, then $\chi(z)\in \U(n)$.

\sm

{\rm c)} A conjugacy class $\frg$ is uniquely determined by the characteristic function
$\chi_\frg(z)$ and the set $\eta_\frg$.   
\end{theorem}

A similar object is a semigroup of conjugacy classes of $Q:=\GL(\infty,\C)$ by a subgroup
$H$ consisting of matrices $
\begin{pmatrix}
 1_\alpha&0\\0& h
\end{pmatrix}
$, for a further discussion see \cite{Dym}, \cite{Haz}, \cite{MH}.

The $\circ$-multiplication is a
 representative of a wide class of operations on sets of conjugacy classes
and double cosets. Such operations  arise in a natural way in representation theory of infinite-dimensional
classical groups, see \cite{Ner-char}, \cite{Ner-faa}, \cite{Ner-haar}. The purpose of this paper is
to give a way for a visualization of such operations. We consider a minimal example but our receipt works
in the wider 
generality of  \cite{Ner-faa}, we explain this in Section 4. 

\sm

{\bf\punct Product of conjugacy classes.%
\label{ss:product-coll}} Fix integers $\alpha\ge 0$ and $m\ge 1$.
Let $N> 0$.
Consider the space $\Mat(\alpha+m N)$. We write its elements as block matrices
\begin{equation}
g=
\begin{pmatrix}
a&b_1&\dots& b_m\\
c_1&d_{11}&\dots &d_{1m}\\
\vdots&\vdots&\ddots&\vdots\\
c_m&d_{m1}&\dots&d_{mm}
\end{pmatrix}
\label{eq:g}
\end{equation}
of size $\alpha+N+\dots+N$.
For an element $u\in\GL(N,\C)$ we denote by
$\iota(u)$ the following matrix
\begin{equation}
\iota(u)=
\begin{pmatrix}
1_\alpha&0&0&\dots& 0\\
0&u&0&\dots &0\\
0&0&u&\dots &0\\
\vdots&\vdots&\vdots&\ddots&\vdots\\
0&0&0&\dots&u
\end{pmatrix}
.
\label{eq:iota}
\end{equation}
Denote by 
$$\cM_N=\cM_N^{\alpha,m}:=\Mat(\alpha+mN)//GL(N,\C)$$ the set of conjugacy classes
of $\Mat(\alpha+mN)$ by $\GL(N,\C)$, i.e.,
\begin{equation}
g\sim \iota(u)g\iota(u)^{-1}
.
\label{eq:sim}
\end{equation}
We call elements of this set by {\it colligations}.

Usually, we  omit superscripts $\alpha$, $m$ from $\cM_N^{\alpha,m}$.

\sm

Next, we  formulate the definition of a colligation
 in another form. Consider the space
$V=\C^m$ and $Z_N=\C^N$. Then an operator
$g$ can be regarded as an operator
$$
\C^{\alpha}\oplus \bigl(V\otimes Z_N)\,\to
\C^{\alpha}\oplus \bigl(V\otimes Z_N)
$$
defined up to a conjugation by an element of 
$\GL(N,\C)$ acting in $Z_N$. 

\sm

We wish to define a canonical  multiplication 
\begin{equation}
\cM_{N_1}\times \cM_{N_2}\to \cM_{N_1+N_2}
\label{eq:circ-1}
\end{equation}
($\alpha$, $m$ are fixed).
First, we consider the case $m=2$. The multiplication
is given by
\begin{multline*}
g\circ h= 
\begin{pmatrix}
a&b_1& b_2\\
c_1&d_{11} &d_{12}\\
c_2&d_{21}&d_{22}
\end{pmatrix}
\circ
\begin{pmatrix}
p&q_1& q_2\\
c_1&d_{11} &d_{12}\\
c_2&d_{21}&d_{22}
\end{pmatrix}
=\\=
\begin{pmatrix}
a&b_1&0& b_2&0\\
c_1&d_{11}&0 &d_{12}&0\\
0&0&1_{N_2}&0&0\\
c_2&d_{21}&0&d_{22}&0\\
0&0&0&0&1_{N_2}
\end{pmatrix}
\begin{pmatrix}
p&0&q_1&0& q_2\\
0&1_{N_2}&0&0&0\\
c_1&0&d_{11}&0 &d_{12}\\
0&0&0&1_{N_2}&0\\
c_2&0&d_{21}&0&d_{22}
\end{pmatrix}
=\\=
\begin{pmatrix}
aa'& |& b_1& ab_1'& & b_1& ab_1'
\\
-& + & -& -&-& - &-
\\
c_1a'&|&  d_{11}& c_1b_1'&& d_{12}& c_1b_2'
\\
c'_1&|&   0 &  d_{11}'&& 0 &  d_{12}'& 
\\
&|&&&&
\\
c_2a'&|&  d_{21}& c_2b_1'&& d_{22}& c_2b_2'
\\
c'_2& |&  0 &  d_{21}'&& 0 &  d_{22}'& 
\end{pmatrix}
\end{multline*}
and we take the conjugacy class containing the latter matrix.

Pass from $m=2$ to an arbitrary $m$ is evident.
However we formulate a
 formal definition, which is valid for all $m$. We extend the  operator $g$
to an operator $\wt g$ in 
$$
\C^\alpha\oplus  \bigl(V\otimes Z_{N_1+N_2}\bigr)=
\Bigl[\C^\alpha\oplus  \bigl(V\otimes Z_{N_1}\bigr)\Bigr]
\oplus \Bigl[V\otimes Z_{N_2}\Bigr]
$$
This operator acts as $g$ on the first summand and as $1$ on the second summand.
Next, we extend $h$ to an operator $\widehat h$ acting in the same way in the space
$$
\C^\alpha\oplus  \bigl(V\otimes Z_{N_1+N_2}\bigr)=
\Bigl[\C^\alpha\oplus  \bigl(V\otimes Z_{N_2}\bigr)\Bigr]
\oplus \Bigl[V\otimes Z_{N_1}\Bigr]
,$$
and set
$$
g\circ h:=
\wt g\cdot\widehat h.
$$
       
\begin{proposition}
The $\circ$-multiplication is
associative, i.e., for any $N_1$, $N_2$, $N_3$
and $g_1\in\cM_{N_1}$, $g_2\in\cM_{N_2}$, $g_3\in\cM_{N_3}$
the following equality holds:
$$
(g_1\circ g_2)\circ g_3= g_1\circ (g_2\circ g_3).
$$
\end{proposition}

This is more or less obvious (see also \cite{Ner-char}).

\sm

{\bf\punct Infinite-dimensional version.%
\label{ss:inf}}
Denote by $\cM_\infty$ the set of conjugacy classes
on $\MAT(\infty)$ with respect to $\GL(\infty,\C)$.
Then $\cM_\infty$ is a semigroup with respect to $\circ$-multiplication.

\sm

{\bf\punct Variants.\label{ss:variants}} 
Denote 
\begin{align*}
\cG\cL_N&=\cG\cL_N^{\alpha,m}=\GL(\alpha+mN,\C)//\GL(N,\C)\\
\cU_N&=\cU_N^{\alpha,m}=\U(\alpha+mN)//\U(N),
\end{align*}
in all cases we consider the equivalence (\ref{eq:sim}), in the unitary
case we suppose that $u$ is unitary. Thus we get the operations
\begin{align*}
\cG\cL_{N_1}^{\alpha,m}\times \cG\cL_{N_2}^{\alpha,m}\,&\to \cG\cL_{N_1+N_2}^{\alpha,m};\\
\cG\cL_{\infty}^{\alpha,m} \times \cG\cL_{\infty}^{\alpha,m} \,&\to \cG\cL_{\infty}^{\alpha,m};\\
\cU_{N_1}^{\alpha,m}\times \cU_{N_2}^{\alpha,m}\,&\to \cU_{N_1+N_2}^{\alpha,m};\\
\cU_{\infty}^{\alpha,m}\times \cU_{\infty}^{\alpha,m}\,&\to \cU_{\infty}^{\alpha,m}.
\end{align*}
In this section $\alpha$, $m$ are fixed and we will omit upper superscripts.

\sm

{\bf\punct The origin of $\circ$-multiplication.%
\label{ss:origin}}
The spaces $\cM_\infty^{\alpha,1}$ and $\cU_\infty^{\alpha,1}$ were discussed above
in Subsection \ref{ss:classic}.
In 80s a family of operations of this type arose in a natural way in representation
theory of infinite-dimensional classical and  infinite symmetric  groups, see \cite{Olsh-symm}, \cite{Olsh-GB},
\cite{Ner-book}, Section
IX.4. For instance, there is a semigroup structure
on double cosets 
$$
\lim\limits_{n\longrightarrow \infty}
\O(n)\setminus \U(\alpha+n)/\O(n).
$$ 
Moreover, this semigroup acts in the space of $\O(\infty)$-fixed vectors 
of unitary representations of $\U(\alpha+\infty)$. 
Big zoo of operations on  $K\setminus G/ K$ (with 'small'
subgroups $K$) arose in \cite{Ner-char}, \cite{Ner-faa}, \cite{Ner-symm},  \cite{GN}, \cite{Ner-umn};
our $\U(\alpha+m\infty)//\U(\infty)$ is inside this zoo.
For another explanation, see \cite{Ner-haar} 
(the $\circ$-product is an infinite-dimensional limit of convolutions of
$\delta$-measures on conjugacy classes).

\sm

{\bf\punct The purpose of the work.%
\label{ss:purpose}}
In \cite{Ner-char} there was obtained an analog of the  characteristic function
for the semigroup $\cU_\infty^{\alpha,m}$ (see also $p$-adic case in \cite{Ner-p}).
The purpose of this work is to obtain spectral data of matrices  visualizing 
the
$\circ$-multiplications and separating points of $\cU_\infty^{\alpha,m}$ 
and of $\cG\cL_\infty^{\alpha,m}$. In fact, for each conjugacy class we assign 
a rational functions $\chi^{(j)}:\Mat(jm)\to \Mat(j\alpha)$ and a divisors $\xi^{(j)}$ in $\Mat(jm)$.
The main statements are Theorems \ref{th:1}, \ref{th:2}.

Notice that the spectral data for collections of matrices were widely discussed,
see different approaches in surveys \cite{Bea}, \cite{Hit}. It seems that our
approach produces another kind of spectral data. On the other hand,
it extends classical approach of operator theory existing  for $m=1$.

\sm

{\bf\punct Structure of the paper.%
\label{ss:structure}}  Section 2 contains 
  a construction of 'spectral data',
and main Theorems \ref{th:1}, \ref{th:2}. Proofs are contained
in Section 3. Some extensions of our construction are discussed in Section 4.

\section{Spectral data}

\COUNTERS

{\bf\punct Categorical quotient.%
\label{ss:quotient-1}} Notice that the spaces $\cM_N$ are non-Hausdorff.
There are many ways to construct Hausdorff spaces from set-theoretical
 quotients. We will use the following approach  (see, e.g., \cite{PV}).

 Let a reductive group $G$ act
on an affine algebraic variety $X$. Consider the algebra $\C[X]^G$ of $G$-invariant regular
functions on $X$. The {\it categorical quotient}
$[X/G]$ is the set of maximal ideals 
of $\C[X]^G$. 

We consider categorical quotients $[\cM_N]$ of $\Mat(\alpha+mN)$
by $\GL(N,\C)$ and $[\cG\cL_N]$ of $\GL(\alpha+mN,\C)$ by $\GL(N,\C)$.

Notice that regular functions on $\Mat(\alpha+mN)$ are polynomials.
The algebra of regular functions on $\GL(\alpha+mN,\C)$ is generated 
by polynomials and $\det(g)^{-1}$.

\begin{proposition}
\label{pr:categorical-1}
{\rm a)} The natural map $[\cM_N]\to [\cM_{N+1}]$ is injective.

\sm

{\rm b)} The natural map $[\cG\cL_N]\to [\cG\cL_{N+1}]$ is injective,
and $[\cG\cL_N]=[\cM_N]\setminus\{g:\,\det g=0\}$.

\sm

{\rm c)} The natural map $\cU_N\to \cU_{N+1}$
is injective.
\end{proposition}

The statement is proved in Subsection \ref{ss:change-N}.

This allows to define the infinite-dimensional 'categorical quotient'
$[\cM_\infty]$ as the inductive limit.
 On the other hand this shows that the topological quotient
$\cU_\infty$ is Hausdorff.

\begin{proposition}
\label{pr:categorical-2}
The $\circ$-multiplication determines a map of categorical
quotients,
$$
[\cM^{\alpha,m}_{N_1}]\times [\cM^{\alpha,m}_{N_2}]\to
[\cM^{\alpha,m}_{N_1+N_2}].
$$
\end{proposition}

The statement is proved in Subsection \ref{ss:quotient}.

\sm

{\bf\punct Characteristic functions.%
\label{ss:scahr-function-main}} Fix $g\in \Mat(\alpha+mN)$ given by (\ref{eq:g}).
Let $S$ range in $\Mat(m)$.
 We write the following equation:
 \begin{equation}
 \begin{pmatrix}
 q\\x_1\\\vdots\\ x_m
 \end{pmatrix}
 =
 \begin{pmatrix}
a&b_1&\dots& b_m\\
c_1&d_{11}&\dots &d_{1m}\\
\vdots&\vdots&\ddots&\vdots\\
c_m&d_{m1}&\dots&d_{mm}
\end{pmatrix}
\begin{pmatrix}
p\\ s_{11}x_1+\dots+s_{1m}x_m\\
\vdots\\
s_{m1}x_1+\dots+s_{mm}x_m
\end{pmatrix}
\label{eq:def-char}
. \end{equation}
Eliminating variables $x_1$, \dots, $x_m$,
we get a relation
$$
q=\chi_g(S)\,p
,$$
where $\chi_g(S)$ is a rational matrix-valued function,
$\Mat(m)\to\Mat(\alpha).
$
We call  $\chi_g$ by {\it the characteristic function} of $g$.

The following statements were obtained in \cite{Ner-char}.

\begin{proposition} 
The characteristic function of $g$ depends only on a $\GL(N,\C)$-conjugacy class
of $g$ {\rm(}and not on $g$ itself{\rm)}.
\end{proposition}

\begin{proposition}
Consider the natural map $I_N:\cM_N\to\cM_{N+1}$.
Then
$$
\chi_{\vphantom{\widehat I}I_N g}(S)=\chi_g(S).
$$ 
\end{proposition}

This is obvious.

\begin{theorem}
$$
\chi_{g\circ h}(S)=\chi_g(S)\chi_h(S).
$$
\end{theorem}

\begin{theorem} 
Let $g$ be unitary. Then

{\rm a)} If $\|S\|\le 1$, then $\|\chi_g(S)\|\le 1$.

\sm

{\rm b)} If $S$ is unitary, then $\chi_g(S)$ is unitary.
\end{theorem}

{\sc Remark.} The space $\B_n$ of $n\times n$ matrices with norm $<1$ is a well-known object, 
it is an Hermitian symmetric space $\U(n,n)/\U(n)\times\U(n)$, also they are called 
homogeneous  Cartan domains
of the first kind. The unitary group $\U(n)$ is a Shilov boundary of this domain. Characteristic functions
are holomorphic maps $\B_m\to \B_\alpha$, which send the Shilov boundary to the Shilov boundary. 
For $m=1$ such functions ('inner functions') are well-known topic of function theory,
see \cite{Garn}, \cite{Pot}. \hfill $\lozenge$

\sm

Denoting $\wt S:=S\otimes 1_N$, we represent
(\ref{eq:def-char}) as
\begin{align*}
q=ap+ b\wt S x\\
x=cp+d\wt S x
\end{align*}
Therefore
$$
x=(1- d \wt S)^{-1}c p
$$
and
\begin{equation}
\chi_g(S)=a+b\wt S(1-d\wt S)^{-1} c
\label{eq:chi-s-tilde}
.
\end{equation}

\sm


{\bf\punct The language of Grassmannians.%
\label{ss:grassman}} Now we reformulate the definition
of characteristic function. For a linear space
$W$ denote by $\Gr(W)$ the set of all subspaces in $W$. For
an even-dimensional space denote by $\Gr^{1/2}(W)$
the set of subspaces having dimension $\frac 12\dim W$.

Consider the Grassmannian
 $\Gr^{1/2}(V\oplus V)$
 in the linear space $V\oplus V$.
For any $L\in \Gr^{1/2}(V\oplus V)$  consider the subspace 
$$L\otimes Z_N\subset (V\oplus V)\otimes Z_N= (V\otimes Z_N)\oplus(V\otimes Z_N) $$
We write the equation
 \begin{equation}
 \begin{pmatrix}
 q\\x_1\\\vdots\\ x_m
 \end{pmatrix}
 =
 \begin{pmatrix}
a&b_1&\dots& b_m\\
c_1&d_{11}&\dots &d_{1m}\\
\vdots&\vdots&\ddots&\vdots\\
c_m&d_{m1}&\dots&d_{mm}
\end{pmatrix}
\begin{pmatrix}
p\\
y_1\\
\vdots\\
y_m
\end{pmatrix}
. 
\label{eq:char-bis}
\end{equation}
Consider the set $\cX_g(L)$ of all $(q,p)$ such that there are
$(x,y)\in L\otimes Z_N$ satisfying (\ref{eq:char-bis}).
The set $\cX_g(L)$ is a subspace in $\C^\alpha\oplus \C^\alpha$.
The following statement is straightforward.

\begin{proposition}
Let $L$ be a graph of an operator $S:V\to V$. Assume that $S$
is a non-singular point of $\chi_g$. Then $\cX_g(L)$
is the graph of the operator $\chi_g(S):\C^\alpha\to\C^\alpha$.
\end{proposition}

Thus $\cX_g$ is a rational map $\Gr^{1/2}(V\oplus V)\to \Gr^{1/2}(\C^\alpha\oplus\C^\alpha)$.

\sm


{\bf\punct The distinguished divisor.%
\label{ss:divisor}}
Generally, the characteristic function does not determine an element of $\cU_N$
(and, also of the space $[\cM_N]$).
Indeed, consider the following  matrix
of size $\alpha+(k+l)+(k+l)$:
$$
\begin{pmatrix}
a&b_1&0& b_2&0\\
c_1&d_{11}&0 &d_{12}&0\\
0&0&e_{11}&0&e_{12}\\
c_2&d_{21}&0&d_{22}&0\\
0&0&e_{21}&0&e_{22}
\end{pmatrix}\in \cU^{\alpha,2}_{k+l}
.
$$
Then its characteristic functions does not depend on the matrix
$\begin{pmatrix}e_{11}&e_{12}\\e_{21}&e_{22} \end{pmatrix}$.
Therefore additional invariants are necessary.

 Let $g\in\Mat(\alpha+mN)$. 
We write the equation 
 \begin{equation}
 \begin{pmatrix}
 x_1\\\vdots\\ x_m
 \end{pmatrix}
 =
 \begin{pmatrix}
d_{11}&\dots &d_{1m}\\
\vdots&\ddots&\vdots\\
d_{m1}&\dots&d_{mm}
\end{pmatrix}
\begin{pmatrix}
s_{11}x_1+\dots +s_{1m} x_m\\
\vdots\\
s_{m1}x_1+\dots +s_{mm} x_m
\end{pmatrix}
\end{equation}
and consider the set $\xi_g$ of all $S$ such that this equation has a nonzero  solution.

Next, we reformulate the definition in a more
precise form.
Consider  the polynomial $p_g$ on $\Mat(m)$ determined by the equation
\begin{multline*}
p_g(S)=
\det(1-d\wt S)=\\=
\left(1_{mN}-
 \begin{pmatrix}
d_{11}&\dots &d_{1m}\\
\vdots&\ddots&\vdots\\
d_{m1}&\dots&d_{mm}
\end{pmatrix}
 \begin{pmatrix}
s_{11}\cdot 1_N&\dots &s_{1m}\cdot 1_N\\
\vdots&\ddots&\vdots\\
s_{m1}\cdot 1_N&\dots&s_{mm}\cdot 1_N
\end{pmatrix}\right)
\end{multline*}


The set $\xi_g$ is the set of zeros of the polynomial
$p_g$. It is more natural to consider 
the set $\xi_g$ as a {\it divisor} (see, e.g. \cite{GH}), i.e. we decompose
$p_g(S)$ as a product of irreducible factors
$$
p_g(S)=\prod_i h_i(S)^{v_i}
,$$ 
where $h_i$ are pairwise distinct. Then we consider the collection of
hypersurfaces%
\footnote{of complex codimension 1.} 
$h_i(S)=0$ with assigned multiplicities $v_i$.

Evidently, 
$$
p_{g\circ h}(S)=p_g(S)\,p_h(S).
$$
Equivalently, 
$$
\xi_{g\circ h}=\xi_g +\xi_h.
$$
The sign '$+$' means that we consider the union of hypersurfaces taking in accounts their multiplicities.

\begin{proposition}
\label{pr:char-det}
The following identity holds
\begin{equation}
\det \chi_g(S)
=\frac{\det\begin{pmatrix}a&-b\wt S\\c&1-d\wt S  \end{pmatrix}}{\det(1-d\wt S)}
\label{eq:char-det}.
\end{equation}
In particular, the denominator of the rational function
$\det \chi_g(S)$ is a divisor of $p_g=\det(1-\wt S)$.
\end{proposition}

{\sc Proof.}  Applying the formula for the determinant of a block matrix 
$$
\det \begin{pmatrix}
      A&B\\C&D
     \end{pmatrix}=\det(D)\det (A-BD^{-1}C),
$$
we get
$$
\det\begin{pmatrix}a&-b\wt S\\c&1-d\wt S  \end{pmatrix}=\det(1-d\wt S)\cdot
\det\Bigl[a+b\wt S(1-d\wt S)^{-1} c\Bigr].
$$
By (\ref{eq:chi-s-tilde}), the expression in the square brackets is the characteristic function.
\hfill $\square$

\sm

Let us pass to the language of Grassmannians.
 Consider another matrix coordinate $\Lambda=S^{-1}$ on Grassmannian.  
  Then the equation for the divisor passes to the form
$$
\det(d-\Lambda)=
\left(
 \begin{pmatrix}
d_{11}&\dots &d_{1m}\\
\vdots&\ddots&\vdots\\
d_{m1}&\dots&d_{mm}
\end{pmatrix}
-
 \begin{pmatrix}
\lambda_{11}\cdot 1_N&\dots &\lambda_{1m}\cdot 1_N\\
\vdots&\ddots&\vdots\\
\lambda_{m1}\cdot 1_N&\dots&\lambda_{mm}\cdot 1_N
\end{pmatrix}\right)=0.
$$ 

Two equations 
$$\det(1-dS)=0,\qquad \det(d-\Lambda)=0$$
 determine a divisor
$\Xi_g$
in $\Gr^{1/2}(V\oplus V)$. Indeed, we have two charts in $\Gr^{1/2}(V\oplus V)$,
one consists of graphs of operators $V\oplus 0\to 0\oplus V$, another
from graphs of operators $ 0\oplus V \to V\oplus 0$. These charts do not
cover the whole Grassmannian, but the complement of their union has codimension 2.
Therefore any hypersurface in $\Gr^{1/2}(V\oplus V)$ has an intersection with 
at least one of the charts. In fact, all hypersurfaces in $\Gr^{1/2}(V\oplus V)$ are observable 
in the chart $S$, except the hypersurface $\det\Lambda=0$
(it is a complement to the chart $S$). 

The degree of the polynomial $p_g$ for a generic $g$ is $mN$.
If the divisor $\Xi_g$ contains the component $\det\Lambda=0$
with multiplicity $l$, then the degree of $p_g$ is $m(N-l)$.

\sm

{\bf\punct Characteristic functions for $\cM_\infty$.%
\label{ss:M-infty}}
We have a natural embedding 
$$
I_N:\cM_N\to\cM_{N+1}
$$
 induced by the embedding
$Z_N\to Z_{N+1}$. Evidently,
$$
\chi_{I_N g}(S)=\chi_{g}(S)
.
$$
Therefore, we can define a characteristic function for elements of $\cM_\infty=\cup_N \cM_N$.
On the other hand, we can repeat the definition of characteristic function
 for $N=\infty$,
this produces the same result.

Next,
$$
p_{I_N g}(S)=p_g(S)\cdot \det (1-S)
.$$
Denoting by $\delta$ the divisor $\det(1-S)=0$ we get
$$
\Xi_{I_Ng}=\Xi_g+\delta
.
$$
Thus we can define a 'divisor' $\Xi_g\subset \Gr^{1/2}(V\oplus V)$ for $g\in\cM_\infty$. 
Notice that its component $\delta$ has multiplicity $\infty$.
Except $\delta$, we have finite number of components of finite multiplicity.

\sm

{\sc Remark.}  In the same way, $\xi_{I_Ng}=\xi_g+\delta$. We can define $\xi_g$
for $g\in\cM_\infty$. However, acting in this way we loss information
concerning the multiplicity of the divisor $\det\Lambda=0$.
\hfill $\square$

\sm

{\bf\punct  Multiple characteristic functions.%
\label{ss:more}}

\begin{proposition}
\label{pr:malo}
 For fixed $\alpha$, $m$ and sufficiently large $N$ the Taylor coefficients of the characteristic function 
 $\chi_\frg$
 and of the polynomial $p_\frg$ at zero do not generate the algebra of $\GL(N)$-invariant polynomials
 on $\Mat(\alpha+mN)$.
\end{proposition}

Proof is given in Subsection \ref{ss:malo}.

\sm

For this reason, we introduce additional invariants.
Let $j=1$, $2$, \dots.
For $g\in\GL(\infty,\C)$ we consider the direct sum
$$g^{[j]}:=g\oplus\dots \oplus g$$
of $j$ copies of $g$. It acts in the space
\begin{multline*}
\Bigl[\C^\alpha\oplus (V\otimes Z_N)\Bigr]\oplus\dots
\oplus
\Bigl[\C^\alpha\oplus (V\otimes Z_N)\Bigr]=
\C^{j\alpha}\oplus \bigl((V \oplus\dots\oplus V)\otimes  Z_N\bigr)
 =\\=
\C^{j\alpha}\oplus \bigl((V \otimes \C^j)\otimes  Z_N\bigr) 
.
\end{multline*}
Thus we get an embedding
$$
\cM_N^{\alpha,m}\to \cM_N^{j\alpha,j m}
.
$$
It is compatible with $\circ$-multiplication.

Now for a given $g\in\cM^{\alpha,m}_N$ we get the collection
of characteristic functions
\begin{equation}
\chi_{g},\, \chi_{g^{[2]}}, \chi_{g^{[3]}},\dots;\qquad \chi_{g^{[j]}}:\Mat(j m)\to\Mat(j\alpha)
\label{eq:data-1}
\end{equation}
and the collection of divisors 
\begin{equation}
\Xi_{g},\, \Xi_{g^{[2]}}, \Xi_{g^{[3]}},\dots; \qquad \Xi_{g^{[j]}}\subset
\Gr\bigl((V\otimes\C^j)\,\oplus\, (V\otimes\C^j)\bigr).
\label{eq:data-2}
\end{equation}

\begin{theorem}
\label{th:1}
{\rm a)} For any $\alpha$, $m$, $N$ the characteristic functions 
{\rm(\ref{eq:data-1})}
and divisors {\rm(\ref{eq:data-2})} uniquely determine an element
of the categorical quotient $[\cM^{\alpha,m}_N]$

\sm

{\rm b)} The characteristic functions 
{\rm(\ref{eq:data-1})}
and divisors {\rm(\ref{eq:data-2})} uniquely determine an element
of  the categorical quotient $[\cG\cL^{\alpha,m}_N]$
\end{theorem}

Proof is contained in Subsections \ref{ss:invariant}--\ref{ss:invariant-gl}.

\begin{theorem}
\label{th:2}
{\rm a)}  The characteristic functions 
{\rm (\ref{eq:data-1})}
and divisors {\rm(\ref{eq:data-2})} uniquely determine an element
of the space $\cU^{\alpha,m}_N$ of conjugacy classes.

\sm

{\rm b)}  The characteristic functions 
{\rm(\ref{eq:data-1})}
and divisors {\rm(\ref{eq:data-2})} uniquely determine an element
of the space $\cU^{\alpha,m}_\infty$ of conjugacy classes.
\end{theorem}

Proof is contained in Subsection \ref{ss:invariant-u}.





\section{Invariants}

\COUNTERS

{\bf\punct Invariants for $\cM^{\alpha,m}_N$.%
\label{ss:invariant}}
Denote
$$
a^{[j]}=a\otimes 1_j, \qquad b^{[j]}=b\otimes 1_j, 
\qquad c^{[j]}=c\otimes 1_j,\quad d^{[j]}=d\otimes 1_j
.$$
Denote by $b_\beta[k]$ the $k$-th row of the matrix $b_\beta$,
by $c_\gamma[l]$ the $l$-th column of the matrix $c_\gamma$.
Let us regard  $\Mat(\alpha+mN)$ as a linear space with
action of $\GL(N,\C)$. A point of the space is  a collection of $m^2$ matrices $d_{ij}$,
$\alpha m$ vectors $c_i[l]$, and $\alpha m$ of covectors $b_j[k]$.
Generators of the algebra of invariants are known (see, e.g., \cite{Pro},
Section 11.8.1). The algebra is generated by the following polynomials
\begin{align}
&\tr d_{\phi_1 \psi_1} d_{\phi_2 \psi_2}\dots d_{\phi_n \psi_n},
\label{eq:invariant-1}
\\
&b_\beta[k] d_{\phi_1 \psi_1} d_{\phi_2 \psi_2}\dots d_{\phi_n \psi_n}   c_{\gamma}[l],
\label{eq:invariant-2}
\\
&a_{\sigma\tau}.
\label{eq:invariant-3}
\end{align}

We wish to show that all the generators can be
expressed in the terms of Taylor coefficients of
$\chi_{g^{[j]}}(S)$ and  $p_{g^{[j]}}(S)$ at zero.

First, consider the expression
\begin{equation}
\ln p_{g^{[j]}}(S)=
\ln
\det(1- d^{[j]}\wt S),
\label{eq:taylor-1}
\end{equation}
where $S$ is an operator in $V\otimes\ \C^j$.

 It is convenient
to think that matrix elements $s_{\phi\psi}^{\mu\nu}$ of $S$ depend on 4 indexes:
indices
$\phi$, $\psi\le m$ are responsible for an operator in the space $\C^m$,
and $\mu$, $\nu\le j$ for operators in the space $\C^j$.
In a neighborhood of $S=0$ we have the following expansion
\begin{multline}
\ln\det (1-d^{[j]}\wt S)=
\sum_{n>0} \frac{(-1)^n}{n} \tr  (d^{[j]}\wt S)^n=\\=
\sum_{\begin{matrix}
n\ge 0,\\ \phi_1,\dots,\phi_n, \psi_1,\dots,\psi_n\le m,
\\
\mu_1,\dots,\mu_n \le j
\end{matrix}}  \frac{(-1)^n}{n}
s_{\phi_1,\psi_1}^{\mu_1,\mu_2} s_{\phi_2,\psi_2}^{\mu_2,\mu_3}\dots s_{\phi_n,\psi_n}^{\mu_n,\mu_1}
\times\\ \times
\tr d_{\psi_1,\phi_2} d_{\psi_2,\phi_3}\dots d_{\psi_n,\phi_1}
.
\end{multline}
Generally, Taylor coefficients (in $s^{\mu\nu}_{\phi\psi}$)  of this series are 
sums of several traces. However, assume that all elements of the  sequence
$$
\mu_1,\dots,\mu_n 
$$
are pairwise distinct. Then the coefficient  at 
$$ s_{\phi_2,\psi_2}^{\mu_2,\mu_3}\dots s_{\phi_n,\psi_n}^{\mu_n,\mu_1}$$
is
$$
\tr  d_{\psi_1,\phi_2} d_{\psi_2,\phi_3}\dots d_{\psi_n,\phi_1}.
$$
We observe that all the invariants (\ref{eq:invariant-1})
are contained in the set of Taylor coefficients of
(\ref{eq:taylor-1}).

Next we consider  the characteristic function
$$
\chi_{g^{[j]}}(S)=a^{[j]}+b^{[j]} \wt S(1-d^{[j]} \wt S)^{-1} c^{[j]}
.
$$
Substituting $S=0$, we get invariants (\ref{eq:invariant-3}).
Next, expand the matrix
\begin{equation}
X:=
b^{[j]} \wt S(1-d^{[j]} \wt S)^{-1} c^{[j]}=\sum_{n=0}^\infty
b^{[j]} \wt S (d^{[j]} \wt S)^n c^{[j]}
\label{eq:X}
\end{equation}
in a Taylor series.
The operator  $X$ acts in $\C^m\otimes \C^j$.  We enumerate its matrix elements
as $x_{kl}^{\gamma\delta}$, where $k$, $l\le \alpha$ and $\gamma$, $\delta\le j$.
The matrix elements are
\begin{multline*}
x_{kl}^{\gamma\delta}=\sum_{\begin{array}{c}n\ge 0\\
                           \phi_1,\dots,\phi_{n+1}\le m\\
                           \psi_1,\dots, \psi_{n+1}\le m\\
                           \mu_2,\dots,\mu_{n}\le j \end{array}}
 s_{\phi_1,\psi_1}^{\gamma,\mu_2} s_{\phi_2,\psi_2}^{\mu_2,\mu_3}\dots s_{\phi_{n+1},\psi_{n+1}}^{\mu_n,\delta}
\times\\ \times
b_{\phi_1}[k] d_{\psi_1,\phi_2} d_{\psi_2,\phi_3}\dots d_{\psi_n,\phi_{n+1}}  c_{\psi_{n+1}}[l]
\end{multline*}
Again, consider a coefficient of the Taylor series at 
$$ s_{\phi_1,\psi_1}^{\gamma,\mu_2} s_{\phi_2,\psi_2}^{\mu_2,\mu_3}\dots s_{\phi_{n+1},\psi_{n+1}}^{\mu_n,\delta}$$
with
 pairwise distinct  $\mu_2$, \dots, $\mu_{n}$ that  do not equal $\gamma$, $\delta$.
The  factors $d_{...}$ in the coefficient and their order are
are uniquely determined.  

This proves Theorem \ref{th:1}.a.


\sm

{\bf\punct Invariants for invertible matrices%
\label{ss:invariant-gl}.}
Now consider the space $\cG\cL_N^{\alpha,m}$.
Invariant regular functions on this space have the form
$p(g)/\det(g)^k$, where $p(g)$ is a polynomial
satisfying
$$
p\bigl(\iota(u)\, g\, \iota(u)^{-1}\bigr)=p(g) ,\qquad\text{where $u\in\GL(N,\C)$},
$$
$\iota(u)$ is given by (\ref{eq:iota}). Thus the algebra of
invariant finctions is generated by the same invariants and $\det(g)^{-1}$.


\sm

{\bf\punct Invariants for unitary group%
\label{ss:invariant-u}.}
Consider invariants of $\U(N)$ on $\U(\alpha+mN)$.
Consider the algebra $\cA$ of functions on $\U(n)$ generated by matrix elements and
$\det(g)^{-1}$. By the Peter-Weyl theorem it is dense in the algebra of continuous 
functions. 

Consider two orbits $\cO_1$, $\cO_2\subset\U(\alpha+mN)$
 of $\U(N)$, consider an invariant continuous function $\phi$ separating
these orbits. Consider a function $\psi\in \cA$ approximating $\phi$,
it satisfies a condition of the form
$$
\psi\Bigr|_{\cO_1}\le a<b\le \psi\Bigr|_{\cO_2}
.
$$
We have
$$
\psi(g)=\frac{p(g)}{(\det g)^k}
,$$
where $p$ is a polynomial.
The average
$$
f(g)=
\int_{\U(N)} \psi\bigl(\iota(h)\,g\,\iota(h)^{-1}\bigr)\,dh
$$
 has the same form
$\frac{q(g)}{(\det g)^k}$ with a polynomial  $q$.
Therefore $f$ has an invariant regular continuation to the 
whole group $\GL(\alpha+mN,\C)$. If $N>\alpha m$,
then all such functions are polynomials 
in Taylor coefficients of $\chi_g(S)$ and $p_g(S)$.


\sm

{\bf \punct Change of $N$.%
\label{ss:change-N}}
Next, we wish to show that the map $I_N:[\cM^{\alpha,m}_N]\to [\cM^{\alpha,m}_{N+1}]$
is an embedding (Proposition \ref{pr:categorical-1}). The map $I_N$ replaces
\begin{align*}
a&\to a;
\\
b_i&\to \begin{pmatrix}b_i&0 \end{pmatrix};
\\
c_l&\to \begin{pmatrix}c_l\\0 \end{pmatrix};
\\
d_{\phi\psi}&\to
 \begin{pmatrix}d_{\phi\psi}&0\\0&0 \end{pmatrix}, \qquad i\ne j;
 \\
 d_{\phi\phi}&\to  \begin{pmatrix}d_{\phi\phi}&0\\0&1 \end{pmatrix}.
\end{align*}
Restrict generators
 (\ref{eq:invariant-1})--(\ref{eq:invariant-3})
 of  invariants on $\Mat(\alpha+m(N+1))$ 
  to $\Mat(\alpha+mN)$. We get the same expressions
 in all the cases except traces of the form
 $$
 \tr d_{\phi_1 \phi_1} d_{\phi_2 \phi_2} \dots d_{\phi_n \phi_n}
, $$
they are shifted by 1. Therefore, restricting a collection
of generators  of $\GL(N+1,\C)$-invariants
of $\Mat(\alpha+m(N+1))$ we get a collection of generators
for algebra of  $\GL(N,\C)$-invariants on $\Mat(\alpha+mN)$.

Thus we proved the desired statement.

\sm


{\bf \punct Categorical quotient.%
\label{ss:quotient}}
Proposition \ref{pr:categorical-2} follows from
Theorem \ref{th:1} (invariants of  $\circ$-products are uniquely determined by invariants
of factors).

\sm

{\bf\punct Proof of Proposition \ref{pr:malo}.%
\label{ss:malo}}
Consider the algebra of $\U(N)$-invariants on $\Mat(\alpha+mN)$.
The number of the standard generators of the degree $l$ has asymptotics
$\sim  \mathrm{const}\cdot(m^2)^l$ as $l\to\infty$.

The Taylor coefficients  of the characteristic function are polynomials of generators.
Number of the Taylor coefficients equals to $\alpha^2 C_{m^2+l-1}^{l-1}$, it is a polynomial in $l$ of 
degree $m^2$.

The minimal (graded) degree of a relation between standard generators (\ref{eq:invariant-1}) is $N+1$,
see \cite{Pro0}, Theorem 4.3.

It remains to choose a sufficiently large $N$ such that number of independent monomials of
the independent
generators (\ref{eq:invariant-1}) of degree $\le N$ is larger than number of monomials of the Taylor
coefficients of degree $\le N$.

\section{Variations}

\COUNTERS

Let $A_1$, \dots, $A_l$ be square matrices. By $\Diag(A_1,\dots, A_n)$
we denote the block-diagonal matrix with blocks $A_1$, \dots, $A_l$.

\sm

{\bf \punct Several matrices determined up to a common conjugation.%
\label{ss:several-matrices}}
Let $G=\U(\infty)$, and $K$ be the stabilizer of the first $\alpha$ basis vectors
as in Subsection \ref{ss:classic}. Consider the product $G^n=G\times\dots\times G$
of $n$ copies of $G$ and
the subgroup $K$ embedded to $G\times\dots\times G$ by diagonal, $h\mapsto (h,\dots,h)$.
Consider the space $G^n//K$ of conjugacy classes 
$$
(g_1,\dots, g_n)\sim (hg_1h^{-1},\dots, hg_nh^{-1}).
$$
 It is equipped with the component-wise product
 \begin{equation}(g_1,\dots, g_n)\ast (r_1,\dots, r_n):=
 (g_1\circ r_1,\dots, g_n\circ r_n),
 \label{eq:tuples-prod}
 \end{equation}
 where  $\circ$-multiplication of matrices is given by the formula
(\ref{eq:classic}).

To be definite, set $n=2$. Consider a pair of matrices
$$\begin{pmatrix} 
  a_1&b_1\\c_1&d_1
 \end{pmatrix},\qquad 
 \begin{pmatrix} 
  a_2&b_2\\c_2&d_2
 \end{pmatrix}
 $$
 determined up to a common conjugation by a matrix $\begin{pmatrix}1&0\\0&u \end{pmatrix}$.
 We compose a matrix 
 $$
 \begin{pmatrix}
  a_1&0&b_1&0\\
  0&a_2&0&b_2\\
  c_1&0&d_1&0\\
  0&c_2&0&d_2
 \end{pmatrix}
 .
 $$
 It is determined up to a conjugation by a matrix
 $$
\Diag(1_\alpha,1_\alpha,u,u)=\Diag(1_{2\alpha},u,u)
.
 $$
 Thus we embedded the semigroup $G^n//K$ to the semigroup $\cU_\infty^{\alpha n,n}$.
 We know that multiple characteristic functions and distinguished divisors
 separate elements of $\cU_\infty^{\alpha n,n}$ and therefore they separate
 elements of $G^n//K$. We get the analog of Theorem \ref{th:2}.
 Analog of theorem \ref{th:1} also holds but it requires a  separate proof.

 \sm
 
 {\bf\punct Another example with conjugacy classes.%
 \label{ss:another}}
 Now consider the group $Q_{\alpha,m,m}$ of finitary unitary matrices of size
 $$
 \alpha+\underbrace{\infty+\dots+\infty}_{\text{$m$ times}}+\underbrace{\infty+\dots+\infty}_{\text{$m$ times}}
 $$
 determined up to conjugations by the group $H$ of unitary matrices of the form
 $$
 \Diag(1,\underbrace{u,\dots,u}_{\text{$m$ times}},
 \underbrace{v,\dots,v}_{\text{$m$ times}})
 ,$$
 i.e., $H\simeq\U(\infty)\times\U(\infty)$.
 The formula for multiplication on $Q_{\alpha,m,m}//H$ is the same as in Subsection \ref{ss:product-coll}.
 
 Now we explain how to obtain a characteristic function. To be definite, set $m=2$.
 We fix two $m\times m$ matrices, $T$, $S$, set
 $$\Lambda=T^{-1},\qquad K=S^{-1},$$
 and write the equation
 \begin{equation}
 \begin{pmatrix}
  q\\x_1\\x_2\\y_1\\y_2
 \end{pmatrix}=
 \begin{pmatrix}
  a_1&b_1&b_2&b_3&b_4\\
  c_1&d_{11}&d_{12}&d_{13}&d_{14}\\
   c_2&d_{21}&d_{22}&d_{23}&d_{24}\\
    c_1&d_{31}&d_{32}&d_{33}&d_{34}\\
     c_1&d_{41}&d_{42}&d_{43}&d_{44}
 \end{pmatrix}
 \begin{pmatrix}
  p\\
  s_{11}x_1+s_{12}x_2\\
  s_{21}x_1+s_{22}x_2\\
  t_{11}y_1+t_{12}y_2\\
   t_{21}y_1+t_{22}y_2
 \end{pmatrix}
 \label{eq:uv}
 \end{equation}
 (big matrix is an element $g\in Q_{\alpha,m,m}$).
 We eliminate variables $s_{ij}$, $t_{ij}$ and get a relation
 $$
 q=\chi_\frg(S,T)\,p.
 $$
 Thus we get a rational function $\Mat(m)\times\Mat(m)\to \Mat(\alpha)$,
 it sends
 $\B_m\times \B_m$ to $\B_\alpha$ and $\U(m)\times\U(m)\to \U(\alpha)$.
 As above, a product of conjugacy classes corresponds to a point-wise product of functions.
 See several proofs  of similar statements in \cite{Ner-char}.
 
 We also define a distinguished divisor 
 $$
 \Xi_g\in \Gr^{1/2}(\C^m\oplus \C^m)\times \Gr^{1/2}(\C^m\oplus \C^m)
 .
 $$
 In the coordinates,  $(T,S)\in\Mat(n)\times\Mat(n)$ it is determined 
 by the equation
 $$
 \det\left[ 1-
  \begin{pmatrix}
  d_{11}&d_{12}&d_{13}&d_{14}\\
   d_{21}&d_{22}&d_{23}&d_{24}\\
    d_{31}&d_{32}&d_{33}&d_{34}\\
     d_{41}&d_{42}&d_{43}&d_{44}
 \end{pmatrix}
 \begin{pmatrix}
  s_{11}&s_{12}&0&0\\
   s_{21}&s_{22}&0&0\\
   0&0&t_{11}&t_{12}\\
   0&0&t_{21}&t_{22}
 \end{pmatrix}
\right]=0
. $$
In coordinates $(\Lambda,K)\in\Mat(m)\times\Mat(m)$ 
we get the equation
 $$
 \det\left[ 
  \begin{pmatrix}
  \lambda_{11}&\lambda_{12}&0&0\\
   \lambda_{21}&\lambda_{22}&0&0\\
   0&0&\kappa_{11}&\kappa_{12}\\
   0&0&\kappa_{21}&\kappa_{22}
 \end{pmatrix}
 -
  \begin{pmatrix}
  d_{11}&d_{12}&d_{13}&d_{14}\\
   d_{21}&d_{22}&d_{23}&d_{24}\\
    d_{31}&d_{32}&d_{33}&d_{34}\\
     d_{41}&d_{42}&d_{43}&d_{44}
 \end{pmatrix}
\right]=0
. $$
These two charts cover the product of Grassmannians up to a subvariety of codimenion
2. This is sufficient to define a divisor.

As in Subsection \ref{ss:more}, we consider elements 
$$
\frg^{[j]}=\underbrace{\frg\oplus \dots\oplus \frg}_{\text{$j$ times}}
\in Q_{j\alpha,jm,jm}// H.
$$

\begin{theorem}
The collection of characteristic functions $\chi_{\frg^{[j]}}$ and the distinguished divisors
$\Xi_{\frg^{[j]}}$ uniquely determine a conjugacy class $\frg\in Q_{\alpha,m,m}//U(\infty)\times\U(\infty)$.
\end{theorem}

A proof repeats considerations of Section 3.

 \sm
 
 {\bf\punct Example with double cosets.%
 \label{ss:double-cosets}}
 Consider the group $G=\U(\infty)$ and its subgroup $L$ consisting of matrices
 $\begin{pmatrix}
   1_\alpha&0\\0&u
  \end{pmatrix}$, where $u$ is real orthogonal. Consider the group $G^n=G\times \dots \times G$
  ($n$ times) and the subgroup $L$ embedded to $G^n$ by diagonal. Consider double cosets
  $L\setminus G/L$,
  i.e., we consider tuples
 $$\left( \begin{pmatrix} 
  a_1&b_1\\c_1&d_1
 \end{pmatrix},\dots, \begin{pmatrix} 
  a_n&b_n\\c_n&d_n
 \end{pmatrix}\right)
 $$
up to the equivalence 
  \begin{multline*}
\left( \begin{pmatrix} 
  a_1&b_1\\c_1&d_1
 \end{pmatrix},\dots, \begin{pmatrix} 
  a_n&b_n\\c_n&d_n
 \end{pmatrix}\right)\simeq
 \\ \simeq
 \left(  \begin{pmatrix}
   1_\alpha&0\\0&u
  \end{pmatrix}\begin{pmatrix} 
  a_1&b_1\\c_1&d_1
 \end{pmatrix} \begin{pmatrix}
   1_\alpha&0\\0&v
  \end{pmatrix},\dots,
  \begin{pmatrix}
   1_\alpha&0\\0&u
  \end{pmatrix}\begin{pmatrix} 
  a_n&b_n\\c_n&d_n
 \end{pmatrix}  \begin{pmatrix}
   1_\alpha&0\\0&v
  \end{pmatrix}\right),
  \end{multline*}
  where $u$, $v\in \O(\infty)$.
   The set $L\setminus G^n/L$ is equipped with component-wise
    multiplication as (\ref{eq:tuples-prod}).
  
  Again, set $n=2$. We define a characteristic function as in \cite{Ner-char}.
  For this purpose, we fix $2\times 2$ matrices $T$, $S$ and write the equation
  {\small
  \begin{equation}
  \begin{pmatrix}
   q_1^+\\
   x_1\\
   q_2^+\\
   x_2\\
   q_1^-\\
   \scriptstyle
   t_{11}x_1+ t_{12}x_2 \\
   q_2^-\\
    \scriptstyle
    t_{21}x_1+ t_{22}x_2 \\
  \end{pmatrix}
  =
  \begin{pmatrix}
 \!\!  \begin{pmatrix} 
  a_1&b_1\\c_1&d_1
 \end{pmatrix}&&&\\
 &  \!\! \!\!\begin{pmatrix} 
  a_2&b_2\\c_2&d_2
  \!\! \end{pmatrix}&&\\ 
   && \!\! \!\! \begin{pmatrix} 
  a_1&b_1\\c_1&d_1
 \end{pmatrix}^{t-1}&\\
 &&& \!\! \!\!\!\! \begin{pmatrix} 
  a_2&b_2\\c_2&d_2
 \end{pmatrix}^{t-1}\\ 
  \end{pmatrix}
    \begin{pmatrix}
   p_1^+\\
    \scriptstyle
   s_{11}y_1+s_{12}y_2\\
   p_2^+\\
    \scriptstyle
    s_{21}y_1+s_{22}y_2\\
   p_1^-\\
   y_1\\
   p_2^-\\
   y_2
  \end{pmatrix}
  \label{eq:big}
  \end{equation}}
  We eliminate variables $t_{ij}$, $s_{ij}$ and get the characteristic function
  $$
   q:= \chi(T,S) p
  ,$$
  where $p$, $q$ are columns,
  $q=(p_1^+,p_2^+,p_1^-,p_2^-)^t,\quad q=(q_1^+,q_2^+,q_1^-,q_2^-)^t$.

Then the product in $L\setminus G/L$ corresponds to the product of characteristic functions
(see \cite{Ner-char}). 

Now we write the equation (\ref{eq:big}) in the form
{\small
\begin{equation}
\begin{pmatrix}
 p_1^+\\
 y_1\\
  p_1^+\\
 y_1\\
 q_1^+\\
 x_1\\
  q_1^+\\
 x_1
\end{pmatrix}
=
\begin{pmatrix}
 \begin{matrix}
  0&0\\0&0
 \end{matrix}&
 \begin{matrix}
  0&0\\0&0
 \end{matrix}&
 {   \begin{pmatrix} 
  a_1&b_1\\c_1&d_1
 \end{pmatrix}^t}
 &
  \begin{matrix}
  0&0\\0&0
 \end{matrix}\\
  \begin{matrix}
  0&0\\0&0
 \end{matrix}&
  \begin{matrix}
  0&0\\0&0
 \end{matrix}&
  \begin{matrix}
  0&0\\0&0
 \end{matrix}&
    \begin{pmatrix} 
  a_2&b_2\\c_2&d_2
 \end{pmatrix}^t\\
    \begin{pmatrix} 
  a_1&b_1\\c_1&d_1
 \end{pmatrix}&
   \begin{matrix}
  0&0\\0&0
 \end{matrix}&
   \begin{matrix}
  0&0\\0&0
 \end{matrix}&
   \begin{matrix}
  0&0\\0&0
 \end{matrix}
 \\
   \begin{matrix}
  0&0\\0&0
 \end{matrix}&
     \begin{pmatrix} 
  a_2&b_2\\c_2&d_2
 \end{pmatrix}&
   \begin{matrix}
  0&0\\0&0
 \end{matrix}&
   \begin{matrix}
  0&0\\0&0
 \end{matrix}
\end{pmatrix}
\begin{pmatrix}
 p_1^-\\
 s_{11} x_1+  s_{12} x_2\\
  p_2^-\\
 s_{21} x_1+  s_{22} x_2\\
  q_1^-\\
 t_{11} y_1+  t_{12} y_2\\
  q_2^-\\
 t_{21} y_1+  t_{22} y_2\\
\end{pmatrix}
\label{eq:big-2}
\end{equation}}
The big matrix is determined up to a conjugation by a matrix of the form
$$
\Diag(1_\alpha,u,1_\alpha,u,1_\alpha,v,1_\alpha,v)
.
$$
Thus we get an embedding of sets 
$$
\iota:L\setminus G^n/L\to Q_{2m\alpha,2m,2m}//H
$$
(this is not a homomorphism of semigroups!).
Moreover,  characteristic functions of $\frg\in L\setminus G^n/L$ and of $\iota(\frg)$ coincide
(compare (\ref{eq:uv}) and (\ref{eq:big-2})). Therefore multiple characteristic functions and multiple distinguished 
divisors
separate double cosets $L\setminus G^n/L$.

{\tt Math.Dept., University of Vienna,

 Nordbergstrasse, 15,
Vienna, Austria

\&

Institute for Theoretical and Experimental Physics,

Bolshaya Cheremushkinskaya, 25, Moscow 117259,
Russia

\&

Mech.Math. Dept., Moscow State University,

Vorob'evy Gory, Moscow

e-mail: neretin(at) mccme.ru

URL:www.mat.univie.ac.at/$\sim$neretin

wwwth.itep.ru/$\sim$neretin}


\begin{thebibliography}{cc}

\bibitem{Bea}
 Beauville,A. {\it Determinantal hypersurfaces}, Michigan Math. J., 48:1 (2000), 39--64.
 
  \bibitem{Bro1}
 Brodskij, V.M.
{\it On operator nodes and their characteristic functions.}
 Sov. Math., Dokl. 12, 696-700 (1971). 
 
 \bibitem{Bro2}
Brodskij, M.S.
{\it Unitary operator colligations and their characteristic functions.}
Russ. Math. Surv. 33, No.4, 159--191 (1978).
 
 \bibitem{Dym}
Dym, H. {\it Linear algebra in action.}
 American Mathematical Society,
 Providence, RI, 2007.
 
 \bibitem{GN}
 Gaifullin A.A.,  Neretin Yu.A.
  {\it Infinite symmetric group and bordisms of pseudomanifolds.}
  Preprint,
    arXiv:1501.04062 

 
\bibitem{Garn}
Garnett, J. B. {\it Bounded analytic functions.} 
  Academic Press, Inc. New York-London, 1981. 

\bibitem{Goh}
Gohberg, I., Goldberg, S., Kaashoek, M. A.
 {\it Classes of linear operators. Vol. II.}
 Birkhauser, Basel, 1993.

\bibitem{GH} 
Griffiths, Ph.; Harris, J.          
{\it Principles of algebraic geometry.} 2nd ed. 
           Wiley Classics Library. New York, NY: John Wiley \& Sons (1994).
           
            \bibitem{Haz}
Hazewinkel, M.           
{\it Lectures on invariants, representations and Lie algebras in systems and control theory.}
   S\'emin. d'alg\`ebre P. Dubreil et M.-P. Malliavin, 35\`eme Ann\'ee, Proc., Paris 1982, 
   Lect. Notes Math. 1029, 1-36 (1983).
       
 
 \bibitem{Hit}
  Hitchin, N. {\it Riemann surfaces and integrable systems},
  in
 {\it  Integrable systems},
   Clarendon Press, Oxford, 1999, 11--52.
   
\bibitem{Ism}
Ismagilov, R.S.,
{\it Elementary spherical functions on the groups $\SL(2,P)$ over
a field $P$, which is not locally compact with respect to
the subgroup of matrices with integral elements.} 
Math. USSR-Izvestiya, 1967, 1:2, 349--380.   

\bibitem{Liv1}
Livshits, M.S.
{\it On a certain class of linear operators in Hilbert space.}
 Mat. Sb., N. Ser. 19(61), 239-262 (1946);
English transl. in
Amer. Math. Soc. Transl. (Ser. 2), Vol. 13, 61-83
 1960.

\bibitem{Liv2}
Livshits, M. S. {\it On spectral decomposition 
of linear nonself-adjoint operators.} 
 Mat. Sbornik N.S.  34(76),  (1954). 145--199.
English transl. in
Amer. Math. Soc. Transl. (Ser. 2), Vol 5, 1957, 67--114.  

\bibitem{MH}
Martin, C.; Hermann, R.           
{\it Applications of algebraic geometry to systems theory: the McMillan degree 
and Kronecker indices of transfer functions as topological and holomorphic system invariants.}
 SIAM J. Control Optimization
             16, 743-755 (1978).
   
   \bibitem{Ner-book}
   Neretin, Yu. A. {\it Categories of symmetries and
 infinite-dimensional groups.}
 Oxford University Press, New York, 1996.
   
   \bibitem{Ner-char} 
Neretin, Yu. A. {\it Multi-operator colligations and multivariate spherical functions.}
Anal. and Math. Physics,  Anal. Math. Phys. 1 (2011), no. 2-3, 121-138.

\bibitem{Ner-faa}
Neretin, Yu. A. {\it
Sphericity and multiplication of double cosets for infinite-dimensional
 classical groups.} Funct. Anal. Appl., 2011,  45:3, 225-239.

\bibitem{Ner-haar}
Neretin, Yu.A.
{\it On degeneration of convolutions of double cosets at infinite-dimensional limit.}
Preprint,  {\tt http://arxiv.org/abs/1211.6149}.

\bibitem{Ner-symm}
Neretin, Yu.
{\it Infinite tri-symmetric group, multiplication of double cosets,
and checker topological field theories.}
Int. Math. Res. Not. IMRN 2012, no. 3, 501–523. 

\bibitem{Ner-p}
Neretin, Yu.A.
{\it Infinite-dimensional $p$-adic groups, semigroups of double cosets,
and inner functions on Bruhat--Tits builldings.} 
Izvestiya: Mathematics, 2015, 79:3, 512--553

\bibitem{Ner-umn}
Neretin, Yu.A.
{\it Infinite symmetric groups and combinatorial constructions of topological field theory type.}
Russ. Math. Surv., 2015, 70:4, 715--773

\bibitem{Olsh-symm}
Olshanski, G. I.
{\it Unitary representations of $(G,K)$-pairs that are connected with the infinite symmetric group 
$S(\infty)$}. 
Leningrad Math. J. 1 (1990), no. 4, 983-1014 

\bibitem{Olsh-GB}   
 Olshanski, G. I. {\it Unitary representations of infinite dimensional pairs and the formalism of R. Howe}
 in {\it Representation of Lie Groups and Related Topics,}  Gordon and Breach, New York, 1990, 269--463.
 
\bibitem{Olsh-CR} 
 Olshanskii, G. I. {\it Caract\`eres g\'en\'eralis\'es du groupe $\U(\infty)$
 et fonctions int\'erieures},|
 C. R. Acad. Sci. Paris. S\`er. 1, 313:1 (1991), 9-12.

\bibitem{PV}
 Popov,  V.L.;  Vinberg E.B.
 {\it Invariant Theory} in {\it Algebraic Geometry IV}, Encyclopaedia of Math. Sci, v. 55, Springer- Verlag, 1994, 137-314.

\bibitem{Pot}
Potapov, V. P. {\it The multiplicative structure of $J$-contractive
 matrix functions.} Trudy Moskov. Mat. Obshchestva. 4 (1955), 125--236;
English transl. 
Amer. Math. Soc. Transl. (2)  15  (1960) 131--243.   

\bibitem{Pro0}
Procesi, C.
{\it The invariant theory of $n\times n$ matrices}.
Advances in Math. 19 (1976), no. 3, 306-381. 

\bibitem{Pro}
Procesi, C.           
{\it Lie groups. An approach through invariants and representations}, NY: Springer, 2007.

\bibitem{FSN}
Sz.-Nagy, B.; Foias, C.          
{\it Harmonic analysis of operators on Hilbert space.}
 Budapest: Akademiai Kiado (1970).


\end{thebibliography}
\end{document}